\documentclass[a4paper,11pt]{amsart}
\usepackage[plainpages=false]{hyperref}
\usepackage{amsfonts,latexsym,rawfonts,amsmath,amssymb,amsthm, mathrsfs, lscape}
\usepackage{verbatim}

\usepackage[all]{xy}
\usepackage{graphicx,psfrag}

\usepackage{array, tabularx}

\usepackage{setspace}

\usepackage{geometry}
\newtheorem{thm}{Theorem}[section]
\newtheorem{cor}[thm]{Corollary}
\newtheorem{lem}[thm]{Lemma}

\newtheorem{prop}[thm]{Proposition}

\newtheorem{conj}[thm]{Conjecture}
\newtheorem{rmk}[thm]{Remark}

\theoremstyle{definition}

\newcommand{\Z}{{\mathbb{Z}}}
\newcommand{\C}{{\mathbb{C}}}
\newcommand{\R}{{\mathbb{R}}}
\newcommand{\Q}{{\mathbb{Q}}}

\newcommand{\p}{{\mathbb{P}}}

\newcommand{\D}{{i \partial \bar{\partial}}}
\newcommand{\Do}{{i \partial_0 \bar{\partial}_0}}

\newcommand{\Dt}{{i \partial_t \bar{\partial}_t}}

\begin{document}
\title[Deformations of nodal KE Del Pezzo surfaces]
{Deformations of nodal K\"ahler-Einstein Del Pezzo surfaces with discrete automorphism groups}

\author{Cristiano Spotti}
\address{IHES, Le Bois-Marie 35, route de Chartres 91440 Bures-sur-Yvette, France}
\email{spotti@ihes.fr / cristiano.spotti@gmail.com}

\date{\today}
\maketitle

\begin{abstract}
In this paper we prove that generic small partial smoothings of K\"ahler-Einstein (KE) Del Pezzo orbifolds with only nodal singularities, and with no non-zero holomorphic vector fields, admit orbifold KE metrics which are close in the Gromov-Hausdorff sense to the original KE metric.
\end{abstract}

\tableofcontents

\section{Introduction}

Let $(X_0,\omega_0)$  be an orbifold KE Del Pezzo surface. In the case $X_0$ is actually a smooth Del Pezzo surface, it is well known that the existence of a K\"ahler-Einstein metric is stable under small deformations of the complex structure. This follows easily from an implicit function argument combined with the fact that all non-holomorphically rigid smooth Del Pezzo surfaces have discrete automorphism group.

However, when  $(X_0,\omega_0)$ has genuine isolated quotient singularities the situation is more subtle: the possibilty of topological changes in the complex deformation families  due to the smoothings of the singularities and the possible existence of non-zero holomorphic vector fields make the study of the stability of the orbifold KE condition  more complicated.

Understanding the stability of the KE condition under smoothings or partial smoothings of the complex structure is important essentially for two reasons. First of all, it would give a clear picture of the KE metric moduli space near its "boundary" and in particular it would provide explicit example of KE metric degenerations. Secondly, the construction of an orbifold KE metric on a partial smoothing direction would give new examples of orbifold KE metrics. Concerning this last point, we should stress that the existence problem for KE metric on orbifold Del Pezzo surfaces is still not fully understood. However, we should note that recent progresses concerning the case of Del Pezzo surfaces with canonical singularities 
have been made \cite{OSS12}. In particular, it contains the case of nodal Del Pezzo surfaces. Thus, even if the existence problem of KE metrics on nodal Del Pezzo surfaces has been solved, we believe that the methods we describe in this paper to study the local behavior of the KE metrics under small deformations  still has some interest in its own, since it provides an essentially explicit description of the Gromov-Hausdorff degeneration process and it is reasonable to think that the method can be used to study higher dimensional situations.

In this paper we are going to investigate the problem of deformations of singular KE Fano varieties in the simplest case, namely the case of KE Del Pezzo orbifolds with nodal (i.e., $A_1$) singularities and discrete automorphism groups. To be more precise, we prove the following result:

\begin{thm}\label{MT} Let  $(X_0,\omega_0)$  be a KE Del Pezzo orbifold. Suppose that
\begin{itemize}
\item All the singularities are nodes (i.e., locally of the form $z_1^2+z_2^2+z_3^2=0)$;
\item $\sharp (Aut(X_0)) < \infty$.
\end{itemize}
Then, if $\pi:\mathcal{X} \rightarrow \Delta \subseteq \C_t$ is a generic (partial) smoothing of $\pi^{-1}(0)=X_0$, $X_t$ admits an (orbifold) KE metric $\omega_t$ for $|t|<<1$. Moreover $(X_t,\omega_t) \rightarrow (X_0,\omega_0)$ in the Gromov-Hausdorff (GH) sense.
\end{thm} 

By the word ``generic'' we mean that if $t$ is the parameter of the base of the smoothing family $\mathcal{X}$ and if $z_1^2+z_2^2+z_3^2=s \in \C^{3}\times \C_s$ is the total family of the versal deformation of the node \cite{KS72}, then $s$ and $t$ are related by $$s=s(t)=Ct+\mathcal{O}(t^2)$$ with $C\neq 0$.

The proof of the Theorem is based on a gluing construction which ``reverses'' the GH degeneration picture. The basic idea consists in smoothing the singularities by replacing a neighborhood of a singular point with a scaled version of the Eguchi-Hanson space. Similar gluing constructions have been considered by many authors, e.g., D. Joyce \cite{J96}, C. Arezzo and F. Pacard \cite{AP06}, S. Donaldson \cite{D10g} and O. Biquard and V. Minerbe \cite{BM11}. However, in our situation the gluing is complicated by the fact that the complex structure changes in our family. This problem is reflected in the  construction of an approximate solution of the K\"ahler-Einstein equation, i.e., a K\"ahler metric $\tilde{\omega}_t \in c_1(X_t)$ close to the original $\omega_0$. Once the existence of a good approximate solution is established, the proof of the existence of a K\"ahler-Einstein metric becomes basically standard.  

The reason why we need the above genericity assumption is essentially technical and it has mainly to do with our choice of function spaces in which to perform the gluing construction. Removing the above genericity assumption gives rise to a concentration of the error term on the variety $X_t$ away from the singularities (more precisely the Ricci potential of the pre-glued metric becomes too big to be controlled). Nevertheless, we believe that the Theorem should be true even without the genericity assumption.

The structure of the paper is as follows. In Section $2$ we construct a K\"ahler metric on a (partial)-smoothing which is Gromov-Hausdorff close to the original KE metric on the central fiber $X_0$. In Section $3$, we deform, using an implicit function argument between appropriate weighted H\"older-Spaces, the constructed approximate solution  to a genuine solution of the Einstein equation, i.e. to a KE metric on $X_t$ (for $|t|<< 1$). Section $4$ is devoted to the analysis of some examples where we can apply our smoothing theorem and finally, in Section $5$, we discuss some possible  generalizations.

\emph{Acknowledgements}. I would like to express my deepest gratitude to my supervisor Professor Simon K. Donaldson for the constant encouragement and the  precious teachings. I would also like to thank Professors Claudio Arezzo and Frank Pacard for the very helpful discussions.

\section{Construction of an approximate solution}

The aim of this section is to construct a  pre-glued metric $\omega_t$, i.e., a metric on $X_t$  GH close to the original metric $\omega_0$, of the form $$\omega_t=\beta_t+\Dt \phi_t \in c_1(X_t),$$ where $\beta_t$ is some background K\"ahler forms on the fibers varying continuously with $t$.

Let $(X_0, \omega_0)$ be a nodal KE Del Pezzo orbifold and let $p$ be a singular point. Take $(\zeta_1,\zeta_2)$ coordinates on $\C^2/\Z_2$ around $p$ so that locally $$\omega_0=\Do (\varphi^1_0)=\Do (|\zeta|^2+\mathcal{O}(|\zeta|^4)).$$ We say that a function $f(\zeta,\bar{\zeta})$ is   $\mathcal{O}(|\zeta|^k)$ if
$ ||\nabla^j f || \leq C |\zeta|^{k-j}$ for $j\geq 0$. 

Identifying as usual  $\C^2/\Z_2= \{z_1^2+z_2^2+z_3^2=0\}$, where $z_i=z_i(\zeta_1,\zeta_2)$ are quadratic expression in $\zeta_i$, it is easy to see that the metric $\omega_0$  can be written as  $$\omega_0=\D (|z|+\mathcal{O}(|z|^2))_{|V_0}$$ in a neighborhood $V_0 \subseteq X_0$ of the origin.

 By  the theory of versal deformations of hypersurface singularities \cite{KS72}, we have that the general deformation of the node is given by  $ z_1^2+z_2^2+z_3^2=t$, with $t \in \C$. Observe that we may assume $t$ to be real (by ``rotating'' the coordinates). 

Thus, using our genericity assumption, we  \emph{identify} a portion of $z_1^2+z_2^2+z_3^2=t \in \R$ (say $|z|\leq C$ for some positive constant C) with a subset $V_t$ of $X_t$ for some real $t$. Of course this argument can be applied to all singular points $p \in Sing(X_0)$.

Moreover, using the relative anticanonical sections,  we can assume that $\mathcal{X}$ embeds in $\p^n\times \Delta$ and the map $\pi$ is given by composing the anticanonical embedding with the projection onto the base. Since  $\mathcal{O}_{\p^n}(1)_{|X_t}=K_{X_t}^{-k}$ (for notational simplicity we assume $k=1$, since the power $k$ does not enter in our argument),  we can easily construct a background K\"ahler metric  ${\beta_t}_{|X_t}:={FS_t}_{|X_t}$ in $c_1(X_t)$.
Observe that the restriction of the Fubini-Study metric on the singular fiber is, in general, not of orbifold type. However, as it will be clear later in the next sections, the explicit behavior of the background metric close to the singularity is not important for our argument.

\subsection{A map between the smoothing and the central fiber}

Assume for simplicity that $Sing(X_0)=\{p\}$. Let  $V_0=\{z_1^2+z_2^2+z_3^2=0\}$ and $ V_{t}=\{w_1^2+w_2^2+w_3^2=t\}\subseteq X_t$ and consider the diffeomorphism 
$$\begin{array}{cccc}
F_t:& V_0\setminus \{|z|^2 \leq \frac{t}{2} \}  & \longrightarrow &V_{t} \setminus \{|w|^2= t\} \\
                                      &  z_i &  \longmapsto & w_i=z_i+\frac{t}{2|z|^2}\bar{z}_i 
\end{array},$$
between  the singular space and its smoothing (away from the singularities and defined on a small region of $X_0$). Observe that topologically, $$  V_0\setminus \{|z|^2 \leq \tfrac{t}{2} \}  \cong V_{t} \setminus \{|w|^2= t\}  \cong \R\times \R \p^3,$$
and that $L_t:= V_{t} \cap \{|w|^2= t\}\cong S^2$. Moreover note that $|z|\mapsto |w|= \sqrt{|z|^2+\frac{t^2}{4|z|^2}}$.

Let us point out that we can assume that $F_t$ extends to a diffeomorphism  $\psi_t:X_0 \setminus \{|z|^2 \leq \frac{t}{2} \}   \longrightarrow  X_t \setminus \{|w|^2= t\}$.

\begin{lem} There exists a  diffeomorphism $$\psi_t:X_0 \setminus \{|z|^2 \leq \tfrac{t}{2} \}   \longrightarrow  X_t \setminus \{|w|^2= t\},$$
which coincides with the above $F_t$ for $|z|^2 \leq 4$. 

Moreover the complex structures $J_t$ of $X_t$ and $J_0$ of $X_0$  satisfy on  $X_0 \setminus \{ |z|^2 \leq 1 \}$, $||\nabla^k_{\omega_0}( \psi^{*}_t J_t - J_0)||_{\omega_0}=\mathcal{O}(t)$ for all $k$ (or with respect to any other metric equivalent to $\omega_0$). 

\end{lem}

\begin{dimo} As we have previously remarked, we may assume that $\mathcal{X} \hookrightarrow \p^n\times \Delta$. Then $\mathcal{X}$ inherits a K\"ahler metric from the metric $\omega_{FS}+\D |t|^2$ on $\p^n \times \Delta$. We define a connection on $\mathcal{X}\setminus Sing(X_0)$ simply by taking the normal directions to the fiber of the map $pr_2:\mathcal{X} \subseteq \p^n\times \Delta \rightarrow  \Delta$. 

Considering the flow given by the lifting of the radial vector field $X:= -\frac{\nabla |t|^2}{|\nabla |t||^2}$, we find a smooth diffeomorphism  from $X_{|t|} \setminus K_{|t|}$ to $X_0 \setminus \{|z|^2 \leq 4\}$, for some small compact subsets of $K_{|t|}$ of $X_{|t|}$. Let $G_t:X_0 \setminus \{|z|^2 \leq 4\} \rightarrow X_t$ be the inverse of the previous map (assume w.l.o.g that $t$ is a point on the real path $[0,1]\subseteq \Delta$).

It follows that on the strip $S:=\{4\leq |z|^2 \leq 5\} \subseteq X_0$ we have defined two diffeomorphisms $F_t$ and $G_t$ onto some region of $X_t$. Let $S^{'} \subseteq S$ be some smaller strip and use the map $F_t$ to identify  $S$ with its image inside $X_t$. Under the above identification, $G_t$ can be seen as a family of diffeomorphisms onto their image
$$\tilde{G}_t: S^{'} \rightarrow S \cong [0,1] \times \R\p^3,$$
which are by construction smoothly isotopic to the identity on $S^{'}$, i.e. $$\tilde{G}_t \rightarrow Id$$ as $t \rightarrow 0$. In particular $\tilde{G}_t$ is given by the flow of a time dependent vector field $v_t$ defined on $S^{'}$, i.e.,
$$\frac{d}{dt}(\tilde{G}_t)= v_t(\tilde{G}_t).$$
Let $\tau: S \rightarrow [0,1]$ be a smooth cut-off function given by 
$$\tau(z)=\chi(|z|^2),$$
where $\chi$ is a smooth increasing function equal to zero for $|z|^2 \leq 4$ and equal to one for $|z|^2 \geq 5$. Let $w_t$ be the vector field defined by $w_t:=\tau v_t$ and $H_t$ its associated diffeomorphism (for small $t$).
By construction $H_t$ is equal to the identity for $|z|^2 \leq 4$ and equal to $\tilde{G_t}$ for $|z|^2 \geq 5$. Thus the diffeomorphism
$$
\psi_t:=  \left\{ \begin{array}{rl}
 F_t\circ H_t & \mbox{ if $|z|^2 \leq 5$}; \\
  G_t &\mbox{ if $|z|^2 \geq 5$},
       \end{array} \right. 
$$
satisfies the desired property. 

The estimate on the complex structures follows by our construction of the diffeomorphism $\psi_t$ which is given by the flow of a smooth vector field. More precisely, we can cover the total family $\mathcal{X}$ ``away'' from a neighborhood of the singularity with a finite number (by compactness) of charts $\mathcal{U}_j$ of the form $(\zeta_1^j,\zeta_2^j,t)$. Since  by construction $\phi_t$ is isotopic to the identity, in the chart $\mathcal{U}_j$ $\psi_t$ is given by $\zeta_i \mapsto \zeta_i +f^j(t,\zeta,\overline{\zeta})$ for a \emph{smooth} $f^j$ satisfying $f^j(0,\zeta,\overline{\zeta})=0$. Then the claim follows simply by considering the Taylor's expansion of $f^j$ in the variable $t$.

 \qed
\end{dimo}

\subsection{Pre-glued metric close to the singularities}

On $V_1:=\{z_1^2+z_2^2+z_3^2=1\}\cong TS^2$ we have the CY conical Eguchi Hanson metric: $$\eta_1:= \D(\sqrt{|z|^2+1})_{|V_1}.$$ Note that $(\sqrt{|z|^2+1} -|z|) \leq C \frac{1}{|z|}$ for $|z|>>1$. Pulling-back by the holomorphic map $w_i=\sqrt{t} z_i$ and scaling the metric, we have a (scaling) family of  RF metrics on $V_t\subseteq  \{w_1^2+w_2^2+w_3^2=t\} \subseteq X_t$,
$$\eta_{\delta,t}=  \frac{\delta^2}{\sqrt{t}} \,\D \left(\sqrt{|w|^2+t}\right)_{|V_t},$$
satisfying $\mbox{Diam}_{\eta_{\delta,t}}(L_t)=\delta$.

Pulling-back at the level of the potential the KE metric $$\omega_0=\D (|z|+\mathcal{O}(|z|^2))_{|V_0}= \D (\varphi^1_0)_{|V_0}$$ using the (inverse) of $\psi_t$, we find a metric $\omega^1_t$ for small $t$, degenerate in a small neighborhood of $L_t$, given by 

$$\omega^1_t:= \D \left( {\psi^{-1}_t}^{*} (\varphi^1_0) \right)_{|V_t}.$$
We show how to glue $\eta_{\delta,t}$ to $\omega^1_t$. Take $\delta^2=\sqrt{t}$ and consider the gluing region $ \delta^\alpha\leq |w|_{|V_t} \leq 2 \delta^\alpha$ (topologically $[0,1] \times \R\p^3)$ for $\alpha \in[0,2)$ (observe that $\{|w|=\delta^2\}\cong L_t$). Define
\begin{itemize}
\item $\varphi_\delta^1(w):=  {\psi^{-1}_t}^{*} (\varphi^1_0)$;
\item $\varphi_\delta^2(w):= \sqrt{|w|^2+\delta^4}.$
\end{itemize}

\begin{lem} On the annulus $ \delta^\alpha\leq |w|_{|V_t} \leq 2 \delta^\alpha$, where $\alpha \in[0,2)$, we have:
 $$|\nabla_{\eta_\delta}^{k}(\varphi^1_\delta-\varphi^2_\delta)|_{\eta_\delta} =\mathcal{O}(\delta^{8-3\alpha-\frac{k\alpha}{2}})+ \mathcal{O}(\delta^{2\alpha-\frac{k\alpha}{2}})+ \mathcal{O}(\delta^{4-\alpha-\frac{k\alpha}{2}}).$$
 The error is minimized for $\alpha=\frac{4}{3}$: $|\nabla_{\eta_\delta}^{k}(\varphi^1_\delta-\varphi^2_\delta)|_{\eta_\delta} =\mathcal{O}(\delta^{\frac{2}{3}(4-k)})$.
\end{lem}

\begin{dimo}
First of all  note that on the annulus region $\eta_\delta \approx \Dt (|w|)_{|V_{\delta^4}}$, as the  metric (i.e., the scaled Eguchi-Hanson metric) is equivalent to the restriction of the (singular at the origin) metric $\D  (|w|)$). In particular this implies that  $|\nabla_{\eta_\delta}^k(|w|^j)|_{{\eta_\delta}} \leq C |w|^{j-\frac{k}{2}}.$

By the definition of $\varphi^1_\delta$ and $\varphi^2_\delta$ and since $|z|\mapsto |w|= \sqrt{|z|^2+\frac{t^2}{4|z|^2}}$ under the diffeomorphism $\psi_t$, it follows that

\begin{itemize}
 \item $|\nabla_{\eta_\delta}^{k}(\varphi^2_\delta- |w|)|_{\eta_\delta} \leq C \delta^4 |w|^{-1-\frac{k}{2}}$;
 \item $|\nabla_{\eta_\delta}^{k}(\varphi^1_\delta- \sqrt{\frac{|w|^2+\sqrt{|w|^4-\delta^8}}{2}} )|_{\eta_\delta} \leq C |w|^{2-{\frac{k}{2}}}$.
\end{itemize}
Thus
$$|\nabla_{\eta_\delta}^{k}(\varphi^1_\delta-\varphi^2_\delta)|_{\eta_\delta}\leq   C \left(|w|^{2-{\frac{k}{2}}} +  \delta^4 |w|^{-1-\frac{k}{2}}\right) + $$ $$+|\nabla_{\eta_\delta}^{k}\left(|w|- \sqrt{\frac{|w|^2+\sqrt{|w|^4-\delta^8}}{2}}\right)|_{\eta_\delta},$$
i.e.,
$$|\nabla_{\eta_\delta}^{k}(\varphi^1_\delta-\varphi^2_\delta)|_{\eta_\delta}\leq   C \left(|w|^{2-{\frac{k}{2}}} +  \delta^4 |w|^{-1-\frac{k}{2}} + \delta^8 |w|^{-3-{\frac{k}{2}}}\right).$$

\qed
\end{dimo}

\begin{lem}Define for $\sqrt{t}=\delta^2$ and $|w|_{|X_t} \leq 2$ $$\tilde{\omega}_{t,\delta}^1:=\Dt \left( \chi_\delta \varphi_\delta^1+(1-\chi_\delta)\varphi_\delta^2 \right),$$
where $\chi_{\delta}:=\chi\left(\delta^{-\frac{4}{3}}|w|\right)$ is a smooth increasing cut-off function supported in $|w| \geq \delta^{\frac{4}{3}}$, identically equal to one for $|w| \geq 2\delta^{\frac{4}{3}}$. Then for $\delta$ (hence $t$) sufficiently small
\begin{itemize}
\item $||\nabla_{\delta}^k (\tilde{\omega}_{t,\delta}^1-\eta_{\delta,t})||_{\eta_{\delta,t}} = \mathcal{O}(\delta^{\frac{4-2k}{3}})$;
\item $\tilde{\omega}_{t,\delta}^1>0$.
\end{itemize}

\end{lem}

\begin{dimo}
It follows immediately from the previous Lemma observing that $||\nabla^k_\delta \chi_\delta||_{\eta_{\delta,t}}=\mathcal{O}(\delta^{-\frac{2k}{3}})$ on the strip $|w| \in [\delta^{\frac{4}{3}},2 \delta^{\frac{4}{3}}]$ and zero otherwise.
\qed\\
\end{dimo}

\subsection{The pre-glued metric away from the singularities and matching}

First of all we take the following  global K\"ahler potential for the singular metric $\omega_0$. Let $s$ be any  non-vanishing local section of $K_{X_0}^{-1}$. Then
$$\varphi_0:=\log \frac{|s|^2_{\beta_0}}{|s|^2_{\omega_0}},$$
is a well defined smooth function on $X_0 \setminus Sing(X_0)$. Observe that $\varphi_0$ is just the (log of) the ratio of the two volume forms.
Later we will need to fix a precise section of $K_{X_0}^{-1}$ near the singularities. The correct choice is to take  the section $\hat \Omega_0$ where $\hat\Omega_0$ is the section which pulls-back using the orbifold chart to $\partial_{\zeta_1} \wedge \partial_{\zeta_2}$ with $(\zeta_1,\zeta_2)$ the local charts where  $\omega_0= \delta_{ij}+\mathcal{O}(|\zeta|^2)$.

Since $\omega_0$ is K\"ahler-Einstein, it is  easy to see that
$$\omega_0=\beta_0+\Do \varphi_0.$$

Now we are ready to define the approximate KE metric away from the singularities (i.e., away from $|w|_{X_t} \leq 1$ on $X_t$):
$$\tilde{\omega}_{t,\delta}^2:=\beta_t +\Dt {\psi_t^{-1}}^{*} \varphi_0>0$$
for $t$ sufficiently small and where $\psi_t$ is the map between the smooth and singular fiber constructed before.

The next goal is to match the metrics $\tilde{\omega}_{t,\delta}^1$ and $\tilde{\omega}_{t,\delta}^2$ and construct a metric
$$\tilde{\omega}_{t,\delta} =\beta_t +\Dt \phi_t \in c_1(X_t).$$ 

Since $\eta_{\delta,t}$ is a CY metric, we know that there exists a (explicit) nowhere vanishing holomorphic $(2,0)$ form $\Omega_t$ such that for  $\sqrt{t}=\delta^2$
$$\eta_{\delta,t}^2= C \Omega_t\wedge\overline{\Omega}_t.$$
Let $\hat{\Omega}_t$ be its dual. Then $\hat{\Omega}_t$ is a trivialization of $K_{X_t}^{-1}$ on $V_t$ which satisfies $|\hat{\Omega}_t|_{\eta_{\delta,t}}=1$ (after normalization).

Define $b_t \in C^{\infty}(V_t)$ to be the function given by 
$$b_t:= |\hat{\Omega}_t|_{\beta_t}^2.$$
Then $\beta_t=-\Dt \log b_t$ on $V_t$. Thus $\tilde{\omega}_{t,\delta}^2=\Dt \left(- \log b_t + {\psi_t^{-1}}^{*} \varphi_0\right)$ on $V_t \cap \{1 \leq |w| \leq 2\}$. Since we have that $\tilde{\omega}_{t,\delta}^1=\Dt \varphi^1_\delta$ on the same region, it is natural to match the two metrics at the level of the potential using a cut-off function $\tau_t:=\tau(|w|_{V_t})$, a smooth decreasing cut-off function supported in $|w| \leq 2 $  identically equal to one for $|w| \leq 1$ (notice that here the strip where the cut-off function is non constant is of ``fixed shape'').

\begin{prop}\label{PG} There exists a pluriharmonic function $p_t$ on $V_t$ such that for $\sqrt{t}=\delta^2$
$$\tilde{\omega}_{t,\delta}:= \left\{ \begin{array}{ll}
\tilde{\omega}_{t,\delta}^1  & |w|_{|V_t} \leq 1 \\
\Dt \left( \tau_t\left(\varphi^1_\delta-p_t\right)+(1-\tau_t) \left(- \log b_t + {\psi_t^{-1}}^{*} \varphi_0\right) \right) &  1\leq |w|_{|V_t} \leq 2\\
\tilde{\omega}_{t,\delta}^2 &  \mbox{otherwise}\\
\end{array}  \right.$$
is a K\"ahler metric in $c_1(X_t)$. More precisely $\tilde{\omega}_{t,\delta}-\beta_t=\Dt \phi_t$ where
$$\phi_t=\tau_t \left( \left( \chi_\delta \varphi_\delta^1+(1-\chi_\delta)\varphi_\delta^2 \right) -p_t +\log b_t \right)+(1-\tau_t) \left( {\psi_t^{-1}}^{*} \varphi_0\right) \in C^\infty(X_t,\R).$$
Moreover we can choose $p_t$ to satisfy $|\nabla_{\eta_\delta}^{k} p_t|_{\eta_{\delta,t}}\leq C |w|^{1-\frac{k}{2}}$. 
\end{prop}

\begin{dimo}

Let $q$ be the singular point of $X_0$. By our hypothesis on the metric $\omega_0$ we know that, in the orbifold chart centered in $q$, $\omega_0$ can be expressed for $|\zeta| \leq \sqrt{2}$ as
$$\omega_0=\Do \varphi^1_0=-\Do \log (|\hat{\Omega}_0|^2_{\omega_0})=Ric (\omega_0),$$
where $\varphi^1_0=|\zeta|^2+\mathcal{O}(|\zeta|^4)$ is the local K\"ahler potential. Thus
$$p_0:= \varphi^1_0+ \log (|\hat{\Omega}_{0}|^2_{\omega_0}),$$
is a $\Z_2$-invariant pluriharmonic real function on  $|\zeta| \leq \sqrt{2}$  vanishing at the origin.
We claim that $p_0= \mathfrak{Re} (h_0)$ where $h_0$ is a $\Z_2$-invariant holomorphic function. Since $d \partial_0 p_0= \bar{\partial}_0 \partial_0 p_0=0$, by the Poincar\'e Lemma there exists a $\Z_2$-invariant holomorphic function $h_0$ vanishing at the origin satisfying
$$d\frac{h_0}{2} =\partial_0 \frac{h_0}{2}=\partial_0 p_0.$$
Then, being $p_0$ real,
$$d  (\mathfrak{Re} (h_0)-p_0)= \partial_0 \frac{h_0}{2} +\bar{\partial}_0 \frac{\bar{h}_0}{2} -\partial_0 p_0 - \bar{\partial}_0 p_0=0.$$
Hence $p_0= \mathfrak{Re} (h_0)$ .

Since $h_0$ is a holomorphic function on $\C^2 / \Z_2$ vanishing at the origin, identifying $\C^2 /\Z_2$ with $w_0^2+w_1^2+w_2^2=0 \in \C^3$, we may assume that $h_0=H(w_1,w_2,w_3)_{|X_0}$ where $H$ is an holomorphic function on $\C^3$ vanishing at the origin (since the node is a normal singularity). We define for sufficiently small $t$:  
$$p_t:= \mathfrak{Re} H_{|V_t},$$
which is a pluriharmonic function on the local smoothing for $|w|_{|V_t} \leq 2$, satisfying the desired estimates.

Now define $a_t:= \varphi^1_\delta-p_t$ and $c_t:=- \log b_t + {\psi_t^{-1}}^{*} \varphi_0$. It is evident that the closed $(1,1)$-form $\tilde{\omega}_{t,\delta}$ is well-defined and that it  would be positive definite as long as 
$$||a_t-c_t||+||d(a_t-c_t)||_{\tilde{\omega}_{t,\delta}^2} \rightarrow 0$$ as $t \rightarrow 0$ on   $1\leq |w|_{|V_t} \leq 2$.

By the definition and the estimates of the diffeomorphism $\psi_t$ onto the (non-collapsing) regions $1\leq |w|_{|V_t} \leq 2$:
$$|a_t-c_t| =   |(\log b_t - {\psi_t^{-1}}^{*}\log b_0)+ ({\psi_t^{-1}}^{*}p_0 -p_t)|=\mathcal{O}(t),$$
(and similarly for the derivatives). Hence $\tilde{\omega}_{t,\delta}>0$.

Finally it follows by a simple computation that the global K\"ahler potential of  $\tilde{\omega}_{t,\delta}$ with respect the background metric $\beta_t$ is given by $\phi_t$ defined in the statement of the Proposition.
\qed \\
\end{dimo}

We should observe that Proposition \ref{PG} has an obvious generalization when we consider smoothings, or partial smoothings, with more than one singularity. In fact, it is sufficient to repeat our argument around all singularities which are smoothed out and note that around points which remain singular we can simply glue the orbifold metric $\omega_0$ to  $\tilde{\omega}_{t,\delta}^2 $ on $1\leq |w|_{V_t} \leq 2$ with a cut-off function $\tau_t$ as we did in the proof of the above Proposition. Note that in this last case a pluriharmonic correction is not needed, since $\varphi_0^1$ can be taken  equal to $-\log |\hat{\Omega}|^2_0$.

Then, for sake of notational simplicity, we continue to argue assuming that $Sing(X_0)=\{p\}$.

Now we define the following function (the Ricci potential) for $\sqrt{t}=\delta^2$:
$$f_{\delta,t}:=\log \left(\frac{|s_t|^2_{\beta_t}}{|s_t|^2_{\tilde{\omega}_{t,\delta}}}\right)-\phi_t,$$
 where $\phi_t$ is the global potential of the metric defined in \ref{PG} and $s_t$ a non-vanishing local section of $K_{X_t}^{-1}$. On $V_t$ we take $s_t$ to be equal to $\hat{\Omega}_t$ previously considered.
Then we have the following Proposition:

\begin{prop}\label{E} The function $f_{\delta,t}$ defined above satisfies:
\begin{itemize}
\item $Ric \, \tilde{\omega}_{\delta,t}=\tilde{\omega}_{t,\delta}+\Dt f_{\delta,t}$;
\item Estimates:

 \begin{itemize}\item $||\nabla^{k}_{\tilde{\omega}_{t,\delta}} f_{\delta,t}||_{\tilde{\omega}_{t,\delta}}= \mathcal{O}(\delta^4)$, on $X_t \setminus \{|w| \leq 1 \cap V_t\}$ (recall $\sqrt{t}=\delta^2)$;
\item   $||\nabla^{k}_{\tilde{\omega}_{t,\delta}} f_{\delta,t}||_{\tilde{\omega}_{t,\delta}}= \mathcal{O}(\delta^{4-2\alpha-\frac{k\alpha}{2}})$, on $\{\delta^{\alpha} \leq |w| \leq 2 \delta^{\alpha}\} \cap V_t$ for $\alpha \in [0,\frac{4}{3}]$;
\item  $||\nabla^{k}_{\tilde{\omega}_{t,\delta}} f_{\delta,t}||_{\tilde{\omega}_{t,\delta}}= \mathcal{O}(\delta^{\alpha-\frac{k\alpha}{2}})$, on $\{\delta^{\alpha}\leq  |w| \leq 2 \delta^{\alpha}\} \cap V_t$ for $\alpha \in [\frac{4}{3},2].$
\end{itemize}
\end{itemize}
\end{prop}

\begin{dimo}
It is immediate to check that $f_{\delta,t}$ defined as above is a Ricci potential, i.e., it satisfies the equation $Ric \, \tilde{\omega}_{t,\delta}=\tilde{\omega}_{t,\delta}+\Dt f_{\delta,t}$.

Now we compute how $f_{\delta,t}$ behaves as we approach the singular fiber:

\begin{itemize}
\item Region $\delta^2 \leq |w|_{|V_t} < \delta^{\frac{4}{3}}$. It follows by definition that
$$f_{\delta,t}=\log b_t -\log|\hat{\Omega}_t|^2_{\eta_{\delta,t}}-\varphi^2_\delta+p_t-\log b_t.$$
Since $\eta_{\delta,t}$ is Ricci Flat, $\log|\hat{\Omega}_t|^2_{\eta_{\delta,t}}=0$. Recalling that $|\nabla_{\eta_\delta}^{k}(\varphi^2_\delta- |w|)|_{\eta_\delta} \leq C \delta^4 |w|^{-1-\frac{k}{2}}$ and $|\nabla_{\eta_\delta}^{k} p_t|_{\eta_{\delta,t}}\leq C |w|^{1-\frac{k}{2}}$, we find that
$$ ||\nabla^{k}_{\tilde{\omega}_{t,\delta}} f_{\delta,t}||_{\tilde{\omega}_{t,\delta}}= \mathcal{O}(\delta^{\alpha-\frac{k\alpha}{2}}),$$
on $\{\delta^{\alpha}\leq  |w| \leq 2 \delta^{\alpha}\} \cap V_t$ for $\alpha \in [\frac{4}{3},2]$.

\item Region $\delta^{\frac{4}{3}} \leq |w|_{|V_t} \leq 2 \delta^{\frac{4}{3}}$ (the glueing region). As before we have that
$$f_{\delta,t}= -\log|\hat{\Omega}_t|^2_{\tilde{\omega}_{\delta,t}}-\chi_\delta \varphi_\delta^1-(1-\chi_\delta)\varphi_\delta^2 +p_t.$$

It follows by  $||\nabla_{\delta}^k (\tilde{\omega}_{t,\delta}^1-\eta_{\delta,t})||_{\eta_{\delta,t}} = \mathcal{O}(\delta^{\frac{4-2k}{3}})$ that $$ \log|\hat{\Omega}_t|^2_{\tilde{\omega}_{\delta,t}}= \mathcal{O}(\delta^{\frac{4}{3}}),$$
Hence $$ ||\nabla^{k}_{\tilde{\omega}_{t,\delta}} f_{\delta,t}||_{\tilde{\omega}_{t,\delta}}= \mathcal{O}(\delta^{\frac{4-2k}{3}}).$$

\item  Region $ \delta^{\frac{4}{3}} < |w|_{|V_t} \leq 1$. We have that
$$f_{\delta,t}= -\log|\hat{\Omega}_t|^2_{\tilde{\omega}_{\delta,t}}-\varphi_\delta^1 +p_t.$$
Recalling that the diffeomorphism between the smoothings and the singular fiber is $\psi_t: w_i  \longmapsto w_i+\frac{t}{2|w|^2} \bar{w}_i$ and the definition of $p_t$, it follows that
$$||\nabla^k_\delta ( -\log|\hat{\Omega}_t|^2_{\tilde{\omega}_{\delta,t}}-\varphi_\delta^1 +p_t)||_{\eta_{\delta,t}}\leq Ct|w|^{-2-\frac{k}{2}}.$$
Thus, since $\sqrt{t}=\delta^2$,  for $|w|=\delta^{\alpha}$ with $0 \leq \alpha < \frac{4}{3}$,
$$||\nabla^{k}_{\tilde{\omega}_{t,\delta}} f_{\delta,t}||_{\tilde{\omega}_{t,\delta}}= \mathcal{O}(\delta^{4-2\alpha-\frac{k\alpha}{2}}).$$

\item Region $X_t \setminus \{ |w|_{|V_t} \leq 1\}$. By definition of the diffeomorphism $\psi_t$, we can assume we have complex structures $J_t$ on $X_t$ and $J_0$ on $X_0$  satisfying on  $X_t \setminus \{ |w|_{|V_t} \leq 1\}$, $||\nabla^k_\delta( J_t -{\psi_t^{-1}}^{*}J_0)||_{\tilde{\omega}_{t,\delta}}=\mathcal{O}(t)$ for all $k$. Then,  since $\sqrt{t}=\delta^2$,
$$||\nabla^{k}_{\tilde{\omega}_{t,\delta}} f_{\delta,t}||_{\tilde{\omega}_{t,\delta}}= \mathcal{O}(\delta^4),$$
for all $k$.

\end{itemize}
\qed \\
\end{dimo}

\section{Deformation to a genuine solution}

Let $\tilde{\omega}_{t,\delta}$ be the approximate KE metric on $X_t$ constructed before, where $t$ and $\delta$ satisfy the relation $\sqrt{t}=\delta^2$. Let us recall the basic properties of this metric needed in this section:
\begin{itemize}
\item $\tilde{\omega}_{t,\delta}=\beta_t+\Dt \phi_t \in c_1(X_t)$ (and $\tilde{\omega}_{t,\delta}$ is of orbifold type if some singularities are kept under partial-smoothing) ;
\item ${\tilde{\omega}}_{t,\delta}$ restricted to $ \{ |w| \leq \delta^{\frac{4}{3}} \} \cap V_t$ is equal to  $\eta_{t,\delta}$ (Eguchi Hanson metric) with vanishing cycle  satisfying  $\mbox{Diam}_{\tilde{\omega}_{t,\delta}}(L_t) \approx \delta$.
\item The Ricci potential $f_\delta$ satisfies:
\begin{itemize}
\item $||\nabla^{k}_{\tilde{\omega}_{t,\delta}} f_\delta||_{\tilde{\omega}_{t,\delta}}= \mathcal{O}(\delta^4)$, on $X_t \setminus \{|w| \leq 1 \cap V_t\}$ (recall $\sqrt{t}=\delta^2)$;
\item   $||\nabla^{k}_{\tilde{\omega}_{t,\delta}} f_\delta||_{\tilde{\omega}_{t,\delta}}= \mathcal{O}(\delta^{4-2\alpha-\frac{k\alpha}{2}})$, on $\{\delta^{\alpha} \leq |w| \leq 2 \delta^{\alpha}\} \cap V_t$ for $\alpha \in [0,\frac{4}{3}]$;
\item  $||\nabla^{k}_{\tilde{\omega}_{t,\delta}} f_\delta||_{\tilde{\omega}_{t,\delta}}= \mathcal{O}(\delta^{\alpha-\frac{k\alpha}{2}})$, on $\{\delta^{\alpha} \leq |w| \leq 2 \delta^{\alpha}\} \cap V_t$ for $\alpha \in [\frac{4}{3},2].$
\end{itemize}
Note that the error is bigger exactly at the gluing region $|w| \approx \delta^{\frac{4}{3}}$, where it behaves like  $\delta^{\frac{4-2k}{3}}$.
\end{itemize}

In order to find a KE metric, we need to solve the following equation:
$$E_{t,\delta}[\varphi]:= \frac{\left(\tilde{\omega}_{t,\delta}+\Dt \varphi\right)^2}{\tilde{\omega}_{t,\delta}^2}-e^{f_\delta -\varphi}=0,$$
where  $\sqrt{t}=\delta^2$ and $\varphi$ is a real valued smooth function on $X_t$. Then $\tilde{\omega}_{t,\delta}+\Dt \varphi>0$ would be the desired KE metric. We claim that it is possible to find a solution, provided that $\delta$  (hence the complex structure parameter $t$) is sufficiently small. 

In order to solve the problem we realize the operator $E_{t,\delta}$ in some weighted-H\"older spaces:

$$E_{t,\delta}: \mathcal{U} \subseteq \mathcal{C}^{2,\gamma}_{\tilde{\omega}_{t,\delta}, \beta}(X_t,\R) \longrightarrow  \mathcal{C}^{0,\gamma}_{\tilde{\omega}_{t,\delta}, \beta-2}(X_t,\R),$$
where $\beta$ is some real number in $(-2,0)$ and $\gamma\in(0,1)$ the H\"older exponent, and $\mathcal{U}$ a suitable neighborhood of the origin.

The Banach space $\mathcal{C}^{k,\gamma}_{\tilde{\omega}_{t,\delta}, \beta}(X_t,\R)$ is for $\beta <0$ defined as follows. As a vector space it is simply $\mathcal{C}^{k,\gamma}(X_t,\R)$. However, instead of using the usual H\"older norm,  we use a weighted norm. First of all we  define a weight function $\rho_t:X_t \rightarrow [\delta,1]$:

\begin{itemize}
\item $\rho_t(|w|_{|X_t})= \delta$, if $|w|_{|X_t}\leq 2 \delta^2$;
\item $\rho_t(|w|_{|X_t})= |w|^{\frac{1}{2}}_{|X_t}$, if $3 \delta^2 \leq|w|_{|X_t}\leq \frac{1}{2}$;
\item $\rho_t(|w|_{|X_t})=1$, if $|w|_{|X_t}\geq 1$;
\item $\rho_t$ is a smoothly increasing interpolation between the above values in the remaining regions. 
\end{itemize}

Then the weighted norm is:

$$||\varphi||_ {\mathcal{C}^{k,\gamma}_{\tilde{\omega}_{t,\delta}, \beta}}:= \sum_{j \leq k} || \rho_{t}^{-(\beta-j)} \nabla_{\tilde{\omega}_{t,\delta}}^j \varphi ||_{L^{\infty}_{\tilde{\omega}_{t,\delta}, \beta}} +  [\varphi]_ {\mathcal{C}^{k,\gamma}_{\tilde{\omega}_{t,\delta}, \beta}},$$
where
$$[\varphi]_ {\mathcal{C}^{k,\gamma}_{\tilde{\omega}_{t,\delta}, \beta}}:= \sup_{p,q | \; d_{t} (p,q) \leq inj_{t} \; p \neq q}\left( \min \{ \rho_t^{-(\beta-j-\gamma)}(p), \rho_t^{-(\beta-j-\gamma)}(q)\}\frac{ ||\nabla_{t}^k \varphi(p) - \nabla_{t}^k \varphi(q)||_{t}}{d_{t}^\gamma(p,q)}\right)$$ 
 (We compute the difference of the two derivatives by  parallel transport along the unique minimal geodesic).

Rewrite $E_{t,\delta}[\phi]$ as:
$$E_{t,\delta}[\varphi]=(1-e^{f_\delta})+\mathcal{D}_{t,\delta}[\varphi]+\mathcal{R}[\varphi].$$
Here $\mathcal{D}_{t,\delta}[\varphi]= \Delta_{\tilde{\omega}_{t,\delta}}\varphi+e^{f_{\delta}}\varphi$ and $\mathcal{R}[\phi]$ contains the non-linearities.

The first input we need is the following (scaled)-Schauder estimate.

\begin{lem} If $\delta$ (hence $t$) is sufficiently small  and $\beta \in (-2,0)$, then  
 $$ ||\varphi||_{\mathcal{C}^{2,\gamma}_{\tilde{\omega}_{t,\delta},\beta}} \leq C \left( ||\varphi||_{{L}^{\infty}_{\tilde{\omega}_{t,\delta}, \beta }}+ ||\mathcal{D}_{t,\delta}[\varphi]||_{\mathcal{C}^{0,\gamma}_{\tilde{\omega}_{t,\delta}, \beta-2}} \right)$$
for all $\varphi \in \mathcal{C}^{2,\gamma}_{\tilde{\omega}_{t,\delta},\beta}(X_t)$ with  a positive constant $C$ independent of $\delta$.
\end{lem}

\begin{dimo}
Take a small $c>0$  such that $B_{\omega_0}(p,c) \subseteq X_0$, with $p\in Sing(X_0)$ is metrically equivalent to the ball of radius $1$ in the flat $\C^2/\Z_2$.

Then, by standard Schauder estimates, we know that the desired estimate holds in the region $\psi_t(X_0 \setminus B_{\omega_0}(p,c)) \subseteq X_t$ for a constant $C=C(c)>0$ fixed. More precisely, we can cover the region with (a finite number) of domains  $D_t^i$ where the geometry is close to the euclidean for all $0\leq t \leq c$. By definition, the weighted H\"older norms are equivalent to the usual H\"older norms. Then the claim follows by the standard Schauder estimates since $\mathcal{D}_{t,\delta}$ is elliptic (with smoothly varying bounded coefficients in $t$). 

It remains to show what happens in the ``collapsing'' region $$X_t \setminus \psi_t(X_0 \setminus B_{\omega_0}(p,c)) \subseteq \{w_0^2+w_1^2+w_2^2=t\}.$$

First of all we pull-back the functions $\varphi_i(w)$ from $w_1^2+w_2^2 + w_3^2=t=\delta^4$ to $z_1^2+z_2^2 +z_3^2=1$ using the map $w_i=z_i \delta^2$. Then define the function
$$\tilde{\varphi}(z):=\delta^{-\beta} \varphi(w(z)),$$
and scale the metric $\tilde{\omega}_{t,\delta}$ so that the diameter of the cycle $L_1$ is equal to a constant  (i.e. consider the metric  $g_{t,\delta}:= \frac{1}{\delta^2} \tilde{\omega}_{t,\delta}$). 

Then our scaled Schauder estimates follow once we prove that
$$ ||\tilde{\varphi}||_{\mathcal{C}^{2,\gamma}_{g_{t,\delta},\beta}(\rho_{g_{t,\delta}} \leq \frac{c}{\delta})} \leq C \left( ||\tilde{\varphi}||_{{L}^{\infty}_{g_{t,\delta}, \beta }(\rho_{g_{t,\delta}} \leq \frac{2c}{\delta})}+ ||\Delta_{g_{t,\delta}}\tilde{\varphi}+\tilde{c}\delta^2\tilde{\varphi} ||_{\mathcal{C}^{0,\gamma}_{g_{t,\delta}, \beta-2} (\rho_{g_{t,\delta}} \leq \frac{2c}{\delta})} \right),$$
 for some constant $C$ independent of $\delta$. Here $c,\tilde{c}$ denote  positive constants and the weighted norm is essentially  given by the sum of seminorms $$[\tilde{\varphi}]_{j,\beta}= \sup || \tilde{\rho}_{g_{t,\delta}}^{-(\beta-j)} \nabla_{g} \tilde{\varphi}||_{g_{t,\delta}},$$
(and similarly for the H\"older seminorm), where $\tilde{\rho}_{g_{t,\delta}}$ is a weight function equal to $1$ on $\{x \in X_t \, | \,\rho_{g_{t,\delta}}(x)\leq 1\}$ ($\rho_{g_{t,\delta}}$ denotes the distance from $L_1$ w.r.t. $g_{t,\delta}$), exactly equal to the distance $\rho_{g_{t,\delta}}$ for  $\{x \in X_t \, |\, \rho_{g_{t,\delta}}(x) \geq 2\}$ and  a smooth increasing interpolating function on the remaining annulus.

On the compact piece  $\{x \in X_t \, | \,\rho_{g_{t,\delta}}(x)\leq R\}$ (for a fixed $R>0$),  the estimate holds by standard Schauder estimates (the Riemannian manifold we are considering is just a compact subset containing $L_1$ of the Eguchi Hanson space. Then we can argue as we did in the region away from the singularities). It remains to show what happens in the regions $R\leq \rho_{g_{t,\delta}}(x) \leq \frac{c}{\delta}$. Of course this region can be covered by balls $B_{g_{\delta},t}(p_i, \frac{r_i}{8})$, where $  R \leq r_i \leq \frac{c}{\delta}$. By scaling the usual Schauder estimates on the euclidean balls of radius $1$ and noting that, since $r_i \leq \frac{c}{\delta}$, the balls $B_{g_{\delta},t}(p_i, \frac{r_i}{8})$ are metrically equivalent to standard balls of radius $\frac{r_i}{8}$ in the euclidean space, we have an estimate of the form
$$|\tilde{\varphi}|_{L^\infty(B_{g_{t,\delta}}(p_i,\frac{r_i}{8}))}+ r_i^{\gamma}[\nabla \tilde{\varphi}]_{C^{0,\gamma}(B_{g_{t,\delta}}(p_i,\frac{r_i}{8}))}+\dots+ r_i^{2+\gamma} [\nabla^2 \tilde{\varphi}]_{C^{0,\gamma}(B_{g_{t,\delta}}(p_i,\frac{r_i}{8}))}$$ 
$$\leq C \left( |\tilde{\varphi}|_{L^\infty(B_{g_{t,\delta}}(p_i,\frac{r_i}{4}))}+ r_i^2 |\Delta_{g_{t,\delta}}\tilde{\varphi}+\tilde{c}\delta^2\tilde{\varphi} |_{L^{\infty}(B_{g_{t,\delta}}(p_i,\frac{r_i}{4}))}+ r_i^{2+\gamma} |\Delta_{g_{t,\delta}}\tilde{\varphi}+\tilde{c}\delta^2\tilde{\varphi} |_{C^{0,\gamma}}  \right)  $$  
for a constant $C$ independent of $t$ (and $\delta$). Then the desired estimate follows by multiplying the above inequality by $r_i^{-\beta}$ and observing that on $B_{g_{t,\delta}}(p_i,\frac{r_i}{4}))$ the weight function is bounded by $\frac{r_i}{2} \leq  \tilde{\rho}_{g_{t,\delta}} \leq 2 r_i$.

\qed \\
\end{dimo}

The following estimate is fundamental:
\begin{prop}\label{I} If $\delta$ (hence $t$) is sufficiently small  and $\beta \in (-2,0)$, then 
$$ ||\varphi||_{\mathcal{C}^{2,\gamma}_{\tilde{\omega}_{t,\delta},\beta}} \leq C \, ||\mathcal{D}_{t,\delta}[\varphi]||_{\mathcal{C}^{0,\gamma}_{\tilde{\omega}_{t,\delta}, \beta-2}} $$
for all $\varphi \in \mathcal{C}^{2,\gamma}_{\tilde{\omega}_{t,\delta},\beta}(X_t)$ with  a positive constant $C$ independent of $\delta$.
\end{prop}

\begin{dimo}
The proof is by contradiction. Assume that the above estimate does not hold. Then there exists a sequence $(\delta_i)$ going to zero and smooth functions $\varphi_i$ on $X_{t_i}$ satisfying $ \infty > ||\varphi_i||_{\mathcal{C}^{2,\gamma}_{i,\beta}}\geq C >0$ and $||\mathcal{D}_{i}[\varphi]||_{\mathcal{C}^{0,\gamma}_{i, \beta-2}}\rightarrow 0$.

Then by the above scaled Schauder estimates we must have that $||\varphi_i||_{L^{\infty}_{\beta}}$ is bounded above from zero. W.l.o.g. we can assume that there is a sequence of points $p_i \in X_{t_i}$ where
$$\rho^{-\beta}_i(p_i)|\varphi_i(p_i)|=1.$$
Now we have three cases depending on how (a subsequence of) $p_i$ converges.

$\mathbf{Case \, 1}$: Assume that $\rho_i(p_i)\geq C> 0$. Then $p_i \rightarrow p_0 \in X_0\setminus \mbox{Sing}(X_0)$. By the upper bound on $||\varphi_i||_{\mathcal{C}^{2,\gamma}_{i,\beta}}$ we can assume that on all compact subsets of $ X_0\setminus \mbox{Sing}(X_0)$ (a subsequence of) $\varphi_i \rightarrow \varphi_\infty$ in $\mathcal{C}^{2,\gamma-\epsilon}$ (Ascoli-Arzela) with $\varphi_\infty(p_0)=c>0$ and $|\rho_0^{-\beta}\varphi_\infty| \leq C$,  and similarly for the derivatives, near the singularity (where $\rho$ denotes the distance from the singularity w.r.t. the KE metric $\omega_0$ on $X_0$).

Moreover we have that by the $C^2$ convergence $$\Delta_0\varphi_\infty+\varphi_\infty=0,$$
on all compact  subsets of $ X_0\setminus \mbox{Sing}(X_0)$  and where $\Delta_0$ denotes the Laplacian w.r.t. the KE metric $\omega_0$. Pulling back locally $\varphi_\infty$ to $\C^2$ using the orbifold covering map and recalling the behavior of $\varphi_\infty$ at the origin, it follows that $\varphi_\infty$ is a weak-solution of $\Delta_0\varphi_\infty+\varphi_\infty=0$ as long as $\beta>-2$. In fact, for all $u \in C^{\infty}_{0}(B(0,R))$
$$\int_{B(0,R) \setminus B(0,\rho)} \varphi_\infty \,(\Delta_0 u + u) dV_{\omega_0}= \int_{B(0,R) \setminus B(0,\rho)} (\Delta_0 \varphi_\infty \, + \varphi_{\infty}) u dV_{\omega_0} + $$  $$+\int_{\partial B(0,\rho)} \varphi_{\infty}\frac{\partial u}{\partial \nu}-u\frac{\partial \varphi_{\infty}}{\partial \nu} d\Sigma, $$ 
where $\frac{\partial}{\partial \nu}$ denotes the normal derivative (at $\partial B(0,\rho)$ w.r.t. the metric $\omega_0$). Then
$$\int_{B(0,R) \setminus B(0,\rho)} \varphi_\infty \,(\Delta_0 u + u) dV_{\omega_0} \leq $$  $$ \leq C_1 \rho^3 \sup_{B(0,\rho)}|\nabla_0 u|\sup_{\partial B(0,\rho)} |\varphi_\infty|+ C_2 \rho^3 \sup_{B(0,\rho)}| u|\sup_{\partial B(0,\rho)} |\nabla_0 \varphi_\infty|\leq $$
$$\leq C\rho^{3}(\rho^{\beta}+\rho^{\beta-1}) \rightarrow 0,$$
as $\rho \rightarrow 0$ if $\beta >-2$. By standard regularity theory of elliptic operators, it follows that the weak solution $\varphi_{\infty}$ is actually (orbifold) smooth.

The above implies by a well-known Bochner identity for KE metric on Fano orbifolds (\cite{T00}) that $(\bar{\partial}\varphi_\infty)^{\sharp}$ must be a holomorphic vector field on $X_0$. Using now the discrete automorphism  hypothesis, we have that $ (\bar{\partial}\varphi_\infty)^{\sharp}=0$. Since $\varphi_{\infty}$ is real valued, it follows that $\varphi_{\infty}$ must be constant, hence identically zero by the equation. However this is in contradiction with  $|\varphi_\infty(p_0)|=c>0$.

Now we investigate the cases when $\rho_i(p_i)\rightarrow 0$.

$\mathbf{Case \, 2}: \frac{\delta_i}{\rho_i(p_i)} \rightarrow C > 0$.

First of all we pull-back the functions $\varphi_i(w)$ from $w_1^2+w_2^2+ w_3^2=t=\delta^4$ to $z_1^2+z_2^2 +z_3^2=1$ using the map $w_i=z_i \delta^2$. Then we define the functions
$$\tilde{\varphi_i}(z):=\delta_i^{-\beta} \varphi_i(w(z)).$$
Moreover, scaling the metric $\tilde{\omega}_{t,\delta}$ so that the diameter of $L_1$ is equal to a constant for all $i$ (i.e. blowing up the metric by $\frac{1}{\delta_i^2}$), we have that:
\begin{itemize}
\item $||\tilde{\varphi_i}||_{C^{2,\gamma}(K)} \leq C$, with respect to the norm induced by the blow up metric, for all compact subsets of $K$ of   $z_1^2+z_2^2 +z_3^2=1$.
\item  $||\tilde{\varphi_i}||_{L^\infty}\leq \frac{C}{1+|z|^{-\frac{\beta}{2}}}$;
\item $|\tilde{\varphi_i}(p_i)|= c>0$ for $p_i$ contained in a compact subset of  $z_1^2+z_2^2 +z_3^2=1$.
\end{itemize}

It  follows by Ascoli-Arzel\'a that (a subsequence of) $\tilde{\varphi_i}\rightarrow \tilde{\varphi}_\infty$ in ${C^{2,\gamma-\epsilon}(K)}$ for all compact subsets, $\lim_{|z|\rightarrow +\infty}   \tilde{\varphi}_\infty(z) = 0$ and $ |\tilde{\varphi}_\infty (p_0)|=c>0$ for some point in $z_1^2+z_2^2 +z_3^2=1$ at finite distance from $L_1$.

Moreover the hypothesis  $||\mathcal{D}_{i}[\varphi_i]||_{\mathcal{C}^{0,\gamma}_{i, \beta-2}}\rightarrow 0$ implies that
$$||\Delta_{\frac{1}{\delta_i^2}\tilde{\omega}_{t,\delta}} \tilde{\varphi_i}+\delta_i^2 e^{\tilde{f}_{\delta_i}}\tilde{\varphi_i}||_{L^\infty(K)}\rightarrow 0$$
for all compact $K$ of $z_1^2+z_2^2 +z_3^2=1$. Hence
$$||\Delta_{\frac{1}{\delta_i^2}\tilde{\omega}_{t,\delta}} \tilde{\varphi_i}||_{L^\infty(K)} \leq ||\Delta_{\frac{1}{\delta_i^2}\tilde{\omega}_{t,\delta}} \tilde{\varphi_i}+\delta_i^2 e^{\tilde{f}_{\delta_i}}\tilde{\varphi_i}||_{L^\infty(K)} +\delta_i^2 e^{\tilde{f}_{\delta_i}}||\tilde{\varphi_i}||_{L^\infty(K)}\rightarrow 0.$$
Since blowing-up $\tilde{\omega}_{t,\delta}$ gives in the limit the Eguchi-Hanson Ricci flat metric $\eta_1$, we obtain
 $$\Delta_{\eta_1}  \tilde{\varphi}_\infty=0,$$
that is $ \tilde{\varphi}_\infty$ is an harmonic function on $z_1^2+z_2^2 +z_3^2=1$ equipped with the Eguchi-Hanson metric. If $\beta<0$, $\lim_{|z| \rightarrow +\infty} \tilde{\varphi}_\infty=0$. Then, by the maximum principle, $\tilde{\varphi}_\infty$ must be identically zero, which is again in contradiction with 
$ |\tilde{\varphi}_\infty (p_0)|=c>0$.

$\mathbf{Case \, 3}: \frac{\delta_i}{\rho_i(p_i)} \rightarrow 0$ (and  $\rho_i(p_i)\rightarrow 0$).

We consider the function $\tilde{\varphi_i}$ defined as  we did in the previous case. By hypothesis it follows that $|\tilde{\varphi_i}|(z(p_i))=c>0$ for a sequence $|z(p_i)|=:R_i^2 \rightarrow +\infty$, where $R_i\delta_i \rightarrow 0$. Blowing-down the metrics by a factor $R_i^{-2}$, the new metrics $\frac{1}{\delta_i R_i} \tilde{\omega}_{t,\delta}$ converge to the flat metric $\eta_0$ on $\C^2 /\Z_2$. Arguing as in the previous section, the functions $\chi_i:= R_i^{-\beta} \tilde{\varphi}_i$ converge (up to a subsequence) to a smooth function $\chi_{\infty}$ on $\C^2 / \Z_2 \setminus \{0\}$ which is harmonic w.r.t. the flat metric. Moreover,
\begin{itemize}
 \item $|\chi_{\infty}(p_0)|=C > 0$, where $p_0$ is a point at distance $1$ from the origin;
  \item $|\chi_{\infty}(p)| \leq C d_{\eta_0}(0,p)^{\beta}$.
\end{itemize}
 The pull-back to $\C^2$ of $\chi_\infty$ is an harmonic function which goes to zero at infinity and which is less singular than the Green's function if $\beta \in(-2,0)$.
Thus it extends smoothly to all $\C^2$. Finally, the maximum principle implies that it must be identically zero. However this is in contradiction with $|\chi_{\infty}(p_0)|=C > 0$.

\qed \\
\end{dimo}

\begin{cor} \label{I2} If $\delta$ (hence $t$) is sufficiently small and $\beta \in (-2,0)$, $$\mathcal{D}_{t,\delta}: \mathcal{C}^{2,\gamma}_{\tilde{\omega}_{t,\delta}, \beta}(X_t,\R) \longrightarrow  \mathcal{C}^{0,\gamma}_{\tilde{\omega}_{t,\delta}, \beta-2}(X_t,\R)$$ is invertible with norm of the inverse independent of $\delta$.
\end{cor}
 
\begin{dimo}Observe that in the non-weighted norm the linearized operator is Fredholm of index zero (it is just Laplacian plus 1). Now for a fixed value of $\delta$ (and $\beta\leq0$), the spaces $ \mathcal{C}^{k,\gamma}_{\tilde{\omega}_{t,\delta}, \beta}(X_t,\R)$ and  $\mathcal{C}^{k,\gamma}_{\tilde{\omega}_{t,\delta}}(X_t,\R)$ are equivalent, i.e. they are the same vector space with equivalent norms. This implies that  $\mathcal{D}_{t,\delta}: \mathcal{C}^{2,\gamma}_{\tilde{\omega}_{t,\delta}, \beta}(X_t,\R) \longrightarrow  \mathcal{C}^{0,\gamma}_{\tilde{\omega}_{t,\delta}, \beta-2}(X_t,\R)$ is also Fredholm of index zero.

By the previous estimate we know that $\mathcal{D}_{t,\delta}$  has no kernel for sufficiently small $\delta$. Thus is also surjective (being of index zero). The fact that the norm of the inverse is uniform is again a consequence of  Proposition \ref{I}. 
\qed\\
\end{dimo}

 To solve the equation $E_{\delta,t}(\varphi)=0$ we use the following version of the Implicit Function Theorem (compare Lemma 1.3 in \cite{BM11}). The proof is elementary. 

\begin{lem}\label{In} Let $E:X \rightarrow Y$ be a differentiable  map between Banach spaces and let $R(x):=E(x)-E(0)-D_0E(x)$ be the non linearities. Assume there exists some positive  constants $L$, $r_0$ and $C$ such that:
\begin{itemize}
\item $|R(x)-R(y)|_Y  \leq L(|x-y|_X)(|x|_X+|y|_X),$ for all $x,y$ in $B_X(0,r_0)$;
\item $D_0E$ is invertible with norm of the inverse bounded by $C$.
\end{itemize}
If  for an $r< \mbox{min}\{r_0,\frac{1}{2LC}\}$ the initial error $|E(0)|_Y \leq \frac{r}{2C}$, then there exists a unique solution of the equation $E(x)=0$ in $B_X(0,r)$.
\end{lem}

In order to apply the above Lemma to our equation, we need to take a closer look at the initial error and at the non linearities.

\begin{lem}
The initial error $E_{\delta,t}(0)=1-e^{f_d}$  is estimated in the weighted norms as 
$$||E_{\delta,t}(0)||_{\mathcal{C}^{0,\gamma}_{\tilde{\omega}_{t,\delta}, \beta-2}(X_t)}=\mathcal{O}(\delta^{\frac{8-2\beta}{3}}),$$
for $\beta \in(-2,0)$.
\end{lem}
\begin{dimo}
We show how to estimate the norm $||E_{\delta,t}(0)||_{L^{\infty}_{\tilde{\omega}_{t,\delta}, \beta-2}(\delta^{\frac{4}{3}}\leq |w|_{|V_t}\leq 1)} $. The other estimates are similar.

 First of all note that at first order the error term is simply given by  the Ricci potential $f_\delta$ (which we know to be small in the point-wise norm by Proposition \ref{E}). Then according to the proof of Proposition \ref{E} and he definition of the weighted norm, we have
$$||E_{\delta,t}(0)||_{L^{\infty}_{\tilde{\omega}_{t,\delta}, \beta-2}(\delta^{\frac{4}{3}}\leq |w|_{|V_t}\leq 1)} \leq C \sup \{\rho_t^{-(\beta-2)}|f_\delta|\} \leq C\delta^4 \sup\{|w|^{-\frac{\beta}{2}-1}\},$$ 
where $\sup$ is w.r.t. points in the region $\delta^{\frac{4}{3}}\leq |w|_{|V_t}\leq 1$. Since $\beta>-2$ the above quantity attains its maximum when $|w| \approx \delta^{\frac{4}{3}}$ (i.e., at the gluing region). Thus $||E_{\delta,t}(0)||_{L^{\infty}_{\tilde{\omega}_{t,\delta}, \beta-2}(\delta^{\frac{4}{3}}\leq |w|_{|V_t}\leq 1)}=\mathcal{O}(\delta^{\frac{8-2\beta}{3}})$ (note that $\delta^4 \leq \delta^{\frac{8-2\beta}{3}}$ for $\delta$ small).  
\qed 
\end{dimo}

By Proposition \ref{I2}, $$||\mathcal{D}_{t,\delta}^{-1}[E_{\delta,t}(0)]||_{\mathcal{C}^{2,\gamma}_{\tilde{\omega}_{t,\delta}, \beta}}=\mathcal{O}(\delta^{\frac{8-2\beta}{3}}).$$
Now  observe that if $||\varphi||_{\mathcal{C}^{2,\gamma}_{\tilde{\omega}_{t,\delta}, \beta}}\leq \mathcal{O}(\delta^{\frac{8-2\beta}{3}})$ then un-weighted norms behave as $||\varphi||_{L^\infty} \leq \mathcal{O}(\delta^{\frac{8+\beta}{3}})$, $||\nabla \varphi||_{L^\infty_\delta} \leq \mathcal{O}(\delta^{\frac{5+\beta}{3}})$ and 
$||\nabla^2 \varphi||_{L^\infty_\delta} \leq \mathcal{O}(\delta^{\frac{2+\beta}{3}})$. In particular these norms go to zero as soon as $\beta > -2$. Thus the  preimage by the linearized operator of the initial error is small in the pointwise norm up to the second derivatives. This is important when we will prove GH convergence of the metrics.

The non linearities are given by the operator
$$\mathcal{R}_{t,\delta} [\varphi]:= \frac{\Dt \varphi \wedge \Dt \varphi}{ \omega_{t,\delta}^2}-e^{f_\delta}\left(\varphi-1+e^{-\varphi}\right).$$

If $||\varphi||_{L^\infty} << 1 $, since $e^{f_d}=\mathcal{O}(1)$ the non-linearities behave as  
$$\mathcal{R}_{t,\delta} [\varphi]=\frac{\Dt \varphi \wedge \Dt \varphi}{ \omega_{t,\delta}^2}-\varphi^2 +\mathcal{O}(\varphi^3),$$
which implies that
$$||\mathcal{R}_{t,\delta}(\varphi_1)-\mathcal{R}_{t,\delta}(\varphi_2)||_{\mathcal{C}^{0,\gamma}_{\tilde{\omega}_{t,\delta}, \beta-2}} \leq C \delta^{(\beta-2)}||\varphi_1-\varphi_2||_{\mathcal{C}^{2,\gamma}_{\tilde{\omega}_{t,\delta}, \beta}}(||\varphi_1||_{\mathcal{C}^{2,\gamma}_{\tilde{\omega}_{t,\delta}, \beta}}+||\varphi_2||_{\mathcal{C}^{2,\gamma}_{\tilde{\omega}_{t,\delta}, \beta}})$$

We are now ready to state and prove the main result:

\begin{prop}\label{M}   If $\delta$ (hence $t$) is sufficiently small and $\beta \in (-2,0)$ then the equation
$$E_{t,\delta}(\varphi_t)=0,$$
admits a (unique) solution  with $||\varphi_t||_{\mathcal{C}^{2,\gamma}_{\tilde{\omega}_{t,\delta}, \beta}} =\mathcal{O}(\delta^{\frac{8-2\beta}{3}})$. 

Moreover $||\nabla^2 \varphi_t||_{L^\infty_\delta} \leq \mathcal{O}(\delta^{\frac{2+\beta}{3}}) \rightarrow 0$ as $\delta \rightarrow 0$.
\end{prop}
\begin{dimo}
We want to apply Lemma \ref{In} to our operators $E_{t,\delta}$. Take $r_0=r_0(\delta)=C_1 \delta^{\frac{8-2\beta}{3}}$. Since $||\varphi||_{L^\infty} \leq   \mathcal{O}(\delta^{\frac{8+\beta}{3}})$, we can use the estimate of the non-linearities $\mathcal{R}_{t,\delta}$ which gives $L=L(\delta)=C_2 \delta^{(\beta-2)}$. In order to apply the Lemma we need that the initial error $E_{t,\delta}(0)$ is much smaller then $\frac{1}{L(\delta)}$, that is
$$\delta^{\frac{8-2\beta}{3}} << \delta^{2-\beta},$$
which is true for $\beta >-2$. Then we can take $r(\delta) \cong r_0(\delta)$ and apply the Lemma. The estimate on the pointwise second derivative follows by the previous observations.

\qed \\
\end{dimo}

Rephrasing the above Proposition, we have that the form
$$\omega_{t,\delta}=\tilde{\omega}_{t,\delta} +\Dt \varphi_t>0$$
is the K\"ahler form of a KE metric provided $\delta$ (hence $t$) is sufficiently small.

Finally we show that the KE metric constructed above converges in the GH topology to the singular metric on the central fiber. Recall the following well-known lemma about GH closeness between compact metrix spaces (compare \cite{BBI01} for a proof):

\begin{lem} \label{GHCC} Let $(X,d_X)$ and $(Y,d_Y)$ be two compact metric spaces.  If $d_{GH}(X,Y) \leq \epsilon$ then there exists a $3\epsilon$-quasi isometry $F:X\rightarrow Y$, i.e.,  a non necessarily continuous map $F:X \longrightarrow Y$  satisfying:
        \begin{itemize}
         \item $|d_X(p,q)-d_Y(F(p),F(q))| \leq 3\epsilon$ for all $p,q$ in $X$
         \item $F(X)$ is $3\epsilon$-dense in $Y$.
        \end{itemize}
Conversely, if there exists a $\epsilon$-quasi isometry $F:X\rightarrow Y$ then $d_{GH}(X,Y) \leq 3\epsilon$.
 
\end{lem}

Using this lemma is then easy to see the following.

\begin{prop} The KE Del Pezzo surface $(X_t,\omega_{t,\delta})$ converges in the Gromov-Hausdorff sense to the original KE Del Pezzo orbifold $(X_0,\omega_0)$.
\end{prop}

\begin{dimo}

It is evident from the construction that the pre-glued metric $\tilde{\omega}_{t,\delta}$ is  GH close to the orbifold metric $\omega_0$. That is $d_{GH}((X_t,\tilde{\omega}_{t,\delta}),(X_0,\omega_0)) \rightarrow 0$ as $\delta$ (hence $t$) goes to zero.

On the other hand it follows from the implicit function argument that
$$||\omega_{t,\delta}-\tilde{\omega}_{t,\delta}||_{\tilde{\omega}_{t,\delta}}=||\Dt\phi_t||_{L^\infty_{\tilde{\omega}_{t,\delta}}}\leq\mathcal{O}(\delta^\frac{2+\beta}{3}),$$
where the norms are computed in the standard point-wise norm w.r.t. the pre-glued metric. In particular, the identity map on $X_t$ is a $C \delta^{\frac{2+\beta}{6}}$-quasi isometry between the KE and the pre-glued metric. 

Then it follows by lemma \ref{GHCC} that
$$d_{GH}((X_t,\omega_{t,\delta}),(X_0,\omega_0)) \leq $$  $$d_{GH}((X_t,\omega_{t,\delta}),(X_t,\tilde{\omega}_{t,\delta}))+ d_{GH}((X_t,\tilde{\omega}_{t,\delta}),(X_0,\omega_0))
\rightarrow  0,$$ as $\delta \rightarrow 0$.

\qed
\end{dimo}

\section{Applications \& Examples}

Before discussing some examples where we can apply the main Theorem \ref{MT}, it is interesting to make some remarks.
First of all, in our considerations on moduli spaces we are principally interested in Del Pezzo orbifolds which appear in GH degenerations of smooth KE Del Pezzo surfaces, that is in smoothable Del Pezzo surfaces. In the next Proposition we see that such smoothable Del Pezzo surfaces with discrete automorphism group cannot have many nodal singularities. 
\begin{prop}Let $X_0$ be a smoothable log Del Pezzo surface with discrete automorphism group and only nodal singularities. Then
$$\sharp\{\mbox{nodes}\}\leq 10-2\, deg(X_0).$$
\end{prop}
\begin{dimo}
Let $\mathcal{X} \rightarrow \Delta_t \subseteq\C$ be a full smoothing of $X_0$ and define 
$$\varphi(t):= \mbox{dim}\,Ext^{0}\left(\Omega_{X_t}^1,\mathcal{O}_{X_t}\right)-\mbox{dim}\,Ext^{1}\left(\Omega_{X_t}^1,\mathcal{O}_{X_t}\right)+\mbox{dim}\,Ext^{2}\left(\Omega_{X_t}^1,\mathcal{O}_{X_t}\right).$$ 
Since $\varphi$ is upper semi-continuous and $X_t$ is a smooth Del Pezzo of degree $deg(X_0)$ for all $t\neq 0$, we find 
 \begin{itemize}
     \item  $\varphi(0) \geq \varphi(t)=\chi (\Theta_{X_t})=2\, deg(X_0)-10$;
      \item $Ext^{0}\left(\Omega_{X_0}^1,\mathcal{O}_{X_0}\right)=0$ since $Aut(X_0)$ is finite and the singularities are normal. 
    \end{itemize}
Thus $$ \mbox{dim}\,Ext^{1}\left(\Omega_{X_0}^1,\mathcal{O}_{X_0}\right) \leq 10-2\, deg(X_0).$$
Finally we can estimate the dimension of the above $Ext^1$ from below by the number of nodes. Using the Grothendieck's local-to-global Ext spectral sequence, the vanishing of $H^{2}(\mathcal{H}om(X,\Omega_{X_0},\mathcal{O}_{X_0}))$, and recalling that $\mathcal{E}xt^{1}\left(\Omega_{X_0}^1,\mathcal{O}_{X_0}\right)$ is a skyscraper sheaf supported on the singular locus with stalk isomorphic to the space of versal deformation of the singularities, we have
$$\begin{array}{rl} \mbox{dim}\,Ext^{1}\left(\Omega_{X_0}^1,\mathcal{O}_{X_0}\right) & \geq  \mbox{dim} \, H^0 \left( \mathcal{E}xt^{1}\left(\Omega_{X_0}^1,\mathcal{O}_{X_0}\right)\right) \\ & =\mbox{dim}\, \left( \oplus_{p \in Sing(X_0)} \C_p  \right) \\ &= \sharp\{\mbox{nodes}\}.\end{array}$$  
\qed
\end{dimo}

\begin{rmk}
 The above Proposition follows also by the \emph{classification} of Del Pezzo surfaces with canonical singularities \cite{BW79}, \cite{D80}, \cite{HW81}, \cite{U83}, \cite{FU86}, \cite{Z88}, \cite{MZ88} and \cite{MZ93}.
\end{rmk}

What can we say about smoothable KE Del Pezzo orbifolds (with the above properties) for each interesting degree $d=(-K_{X_0})^2$?  It is known that there are smoothable del Pezzo quartics, i.e., degree equal to $4$, with no holomorphic vector fields and with only one or two nodal singularities. On the other hand,  it follows by the Mabuchi-Mukai result \cite{MM90} that all KE Del Pezzo surfaces appearing in the GH compactification must have holomorphic vector fields. Even if in this case we cannot use our Theorem to study the local behavior of the moduli space at the boundary points, we can combine it with the Mabuchi-Mukai result to prove the following observation:

\begin{cor} \label{NKE} Let $X$ be a Del Pezzo quartic with only  a nodal singularity and with discrete automorphism group. Then $X$ does not admit a KE metric.
\end{cor} 
   
\begin{dimo} 
 By the classification of intersections of two quadrics (compare for example \cite{AL00}), we can assume that $X_0$ is of the form
\begin{displaymath}
\left\{ \begin{array}{ll}
2 x_0 x_1+x_2^2 +x_3^2+x_4^2=0; & \\
 2 x_0 x_1 +x_0^2 +  x_2^2+a x_3^2+b x_4^2=0, & \end{array} \right.
\end{displaymath}
for generic $a,b\in \C$. Observe that $Aut(X_0)$ is finite. In order to see this, note that the above $X_0$ can be degenerated, by the one-parameter subgroup $T_s$, given by $x_0 \rightarrow  s x_0$, $x_1 \rightarrow s^{-1} x_1$,  to a KE Del Pezzo quartic $X_c$ with two nodes, where the two defining equations are both diagonalized (after an obvious change of coordinates). The automorphism group of $X_c$ can be easily computed by looking at the dimension of the Lie group of matrix in $SL(5,\C)$ which fix the defining equations (note that this group is the full automorphism group since $-K_{X_C}$ is very ample). By a simple computation, it follows that the Lie algebra is one dimensional. Hence $\mbox{dim}_{\C}Aut(X_0)< \mbox{dim}_{\C}Aut(X_c)=1$.   

A generic deformation (smoothing) $X_t$ is given by
\begin{displaymath}
\left\{ \begin{array}{ll}
2x_0 x_1+x_2^2 +x_3^2+x_4^2=0; & \\
2x_0 x_1 +x_0^2 +t x_1^2+  x_2^2+a x_3^2+b x_4^2=0. & \end{array} \right.
\end{displaymath}
If $X_0$  has a KE metric, then by Theorem \ref{MT}  $(X_t,\omega_t)$  degenerates in the GH sense to $X_0$. In particular $X_0$ must appear in the GH compactification. However this is forbidden by the Mabuchi-Mukai result \cite{MM90}.

\qed
\end{dimo}

As we pointed out in the proof of the previous Corollary, an $X_0$ satisfying the hypothesis of the Corollary can be degenerated by a one parameter subgroup $T_s$ to a KE Del Pezzo quartic $X_c$ with two nodes (and with $\C^{*} \subseteq Aut(X_c)$). In particular the Donaldson-Futaki invariant $DF(T_s.X_0)=0$ (since the central fiber admits a KE metric) and $T_s.X_0$ is  a non-trivial test configuration.
We should emphasize  that Stoppa argument for proving that the existence of a KE (cscK) metric implies K-polystability \cite{S09} doesn't extend easily in the case of singular spaces. However, in \cite{B12} the author showed that ``KE implies K-polystability'' is indeed true for $\Q$-Fano varieties. Thus Corollary \ref{NKE} follows also by this more general result.
 
It is more interesting to see what happens for degree three Del Pezzos, i.e., cubic surfaces.  It follows by the well-known classical GIT pictures for cubic surfaces that cubics with nodal singularities  have no holomorphic vector fields. The (unique) cubic with the maximal number of nodal singularities is the so-called Cayley's cubic given by the equation $$C_0:=\{xyz+yzw+zwx+wxy=0\} \subseteq \p^3.$$
Moreover, it is known that  that $C_0$ admits an orbifold KE metric (e.g., \cite{C09}, or simply by noting that it is a quotient of the KE Del Pezzo surface of degree $6$). Thus, considering generic deformations of $C_0$ of the form
$$C_t:=\{xyz+yzw+zwx+wxy+tx^3+(ty^3)=0\} \subseteq \p^3,$$ we can use our Theorem to prove the following:
\begin{cor}\label{ex}Some cubic surfaces with two or three nodal singularities admit KE orbifold metrics. 
\end{cor}

Finally, we also recall that there are several examples of Del Pezzo surfaces of degree $1$ and $2$ to which our Theorem applies.

\section{Further Discussions}

Our main Theorem \ref{MT} admits some natural generalizations. In particular, it is  interesting to analyze what happens in the following three situations:
more general two dimensional quotient singularities, presence of continuous family of automorphisms and higher dimensional (nodal) Fano varieties.

\subsection{More general quotient singularities}

Let $X_0$ be a degree $d$ KE orbifold with $\Q$-Gorenstein smoothable singularities. Recall that an algebraic variety admits a $\Q$-Gorenstein smoothing if it is a central fiber in an (analytic) flat deformation $\mathcal{X} \rightarrow \Delta$ over the disc $\Delta \subseteq \C^*$ such that $-K_{\mathcal{X}}$ is $\Q$-Cartier. In dimension two, singularities which admit (local) $\Q$-Gorenstein smoothings are classified \cite{KSB88}. These singularities are known in the literature as $T$-singularities. Beside canonical singularities, i.e., quotients of $\C^2$ by finite subgroup of $SU(2)$ acting freely away from the origin, $T$-singularities are particular cyclic quotients of $A_k$ canonical singularities.

From the global prospective, it is known that there are no local-to-global obstructions to smooth of Del Pezzo surfaces with $T$-singularities \cite{HP10}, that is, every local smoothing of a $T$-singularities can be realized in a global $\Q$-Gorenstein deformation of $X_0$. 

It is known that (local) smoothing of $T$-singularities admit Asymptotically Conical Calabi Yau metrics. In the case of canonical singularities, this follows by the seminal paper of P. Kronheimer \cite{K89}. More recently, it has been shown by Y. Suvaina in \cite{S12} that Kronheimer metrics descend to the finite quotients which gives $\Q$-Gorenstein smoothings of $T$-singularities.

By the above discussion, it follows that there should be ways to glue  these  Asymptotically Conical Calabi Yau metrics to a KE orbifold with $\Q$-Gorenstein singularities in order to find nearby KE metrics on smoothings of $X_0$.

Finally, we should remark that the above discussion can be also generalized to polarized $\Q$-Gorenstein deformations of constant scalar curvature K\"ahler metrics \cite{S12}\footnote{During the preparation of this article, the author has been informed that O. Biquard and Y. Rollin are studying smoothings of cscK orbifold metrics, in the case of discrete automorphisms, in full generality.}.

\subsection{Non trivial holomorphic vector fields}

If a KE orbifold Del Pezzo $X_0$ with $T$-singularities admits holomorphic vector fields, the problem of finding nearby KE metrics on its (partial)-smoothings is in general obstructed. An easy example of this phenomenon can be seen in the Corollary \ref{NKE}, where it is shown that intersections of two quadrics with only one  nodal singularity, which in particular are small partial smoothings of KE orbifolds with $\C^\ast \subseteq Aut_0(X_0)$, do not admit KE metrics. This example is a special case of the general result: not $K$ polystable $\Q$-Fano varieties does not admit a KE metric. In order to deal with these obstructions, we believe that for small enough partial smoothing of KE Del Pezzo orbifolds with $T$-singularities the problem can be completely understand in term of ``local'' geometry, i.e., by studying the action of the automorphism group of $X_0$ on the space of small $\Q$-Gorenstein partial smoothings. In the case of smooth KE manifolds (more generally cscK), this has been proved by G.Sz\'ekelyhidi in \cite{S10}.
In order to make our discussion a bit more clear, we state now a conjecture on what this expected picture should be.  To make the conjecture more precise, we need to use some deformation theory of singular analytic spaces. The non expert reader can think to $\mbox{E}xt^1(\Omega^1_{X_0}, \mathcal{O}_{X_0})$ as the vector space parameterizing infinitesimal deformations (in the smooth case it corresponds to $H^1(X_0, TX_0)$ used in  Kodaira-Spencer theory). For more information the reader may consult \cite{SP12}, and the references therein.

\begin{conj}
 Let $(X_0,\omega_0)$ be a KE orbifold with $T$-singularities. Denote with $Kur(X_0) \subseteq TDef(X_0)$  the versal Kuranishi space of $\Q$-Gorenstein deformations of $X_0$, where $TDef(X_0)$ is  a vector subspace of $\mbox{E}xt^1(\Omega^1_{X_0}, \mathcal{O}_{X_0})$ parameterizing infinitesimal $\Q$-Gorenstein deformations, and consider the action of the reductive group $Aut_0(X_0)$ on the vector space $TDef(X_0)$. Let $v\in Kur(X_0)$. Then for $|v|$ sufficiently small we have
\begin{itemize}
 \item if $v$ is GIT polystable (i.e., the orbit is closed) for the action of $Aut_0(X_0)$ on $TDef(X_0)$ then $X_v$ admits a KE (orbifold) metrics.
 \item if  $v$ is not GIT polystable, then $X_v$ does not admits a KE (orbifold) metric. 
\end{itemize}

\end{conj}

 We remark that from an analytic point of view the obstruction is given by the non uniform invertibility of the linearized Monge-Amp\`ere operator (basically the operator given by the Laplacian plus one). Observe that on $X_0$ the Kernel (and the coKernel) of the linearized Monge-Amp\'ere equation at the orbifold KE metric is exactly given by the Lie algebra of orbifold holomorphic vector fields, compare for example \cite{T00}. Thus in order to address the above conjecture it will be essential to understand how this analytic obstruction is indeed related to the algebraic stability conditions conjectured to be equivalent to the existence of a KE orbifold metric.

The above discussion have of course a natural generalization to the case of polarized partial smoothings of $\Q$-Gorenstein smoothable cscK orbifolds with non trivial holomorphic vector fields. However, recall that in general in this case there are local-to-global obstruction (i.e., $H^2(X_0,Hom(\Omega^1_{X_0}, \mathcal{O}_{X_0}))$ does not always vanish) and $Aut_0(X_0)$ is not reductive (although one can always restrict to the (reductive) complexification of the group of isometries).

\subsection{Higher dimensions}

Another interesting situation to study is the case of higher dimensional KE $(X_0,\omega_0)$ with nodal singularities, i.e., singularity locally analytically equivalent to $\sum_{i} z_{i}^2=0$. As it should be clear by a careful inspection of this paper, our method for constructing a nearby KE metric on a (partial) smoothing of $X_0$ does not heavily depend on the two dimensional hypothesis, provided we have ``good'' local models to glue in and a careful comprehension of the asymptotics of the metrics involved in the gluing construction.

The asymptotic conical CY metrics one should consider in this case are the so-called Stenzel metrics on the smoothing of the nodes \cite{S93}. In order to perform a gluing construction is then fundamental to understand the behavior of the KE metric (a-priori to be consider only a weak-solution in the sense of pluripotential theory \cite{BBEGZ11}) on the singular space $X_0$ at the singularities. The expected picture is that the metric $\omega_0$ is asymptotic, at some rate, to the well-known Sasaki-Einstein metrics on the node (i.e. $\D (|z|^{2(1-\frac{1}{n})})_{|\{\sum_{i} z_{i}^2=0\}} $. However, at present, there are not known examples where such asymptotic behavior has been established. On the other hand, we should recall that much more it is known about the asymptotic of the Stenzel metric at infinity \cite{CH12}.

In conclusion, under correct assumptions on the decay rates of the metrics, it is reasonable to believe that the method used for the two dimensional situation carries over in the more general higher dimensional setting. This will be the object of future investigations.

\bibliographystyle{plain}

\bibliography{KEnodalDP} 
\end{document}